\documentclass%%[draft]
{article}
\usepackage[T2A]{fontenc}
\usepackage[cp1251]{inputenc}
\usepackage[english,russian]{babel}
\usepackage[tbtags]{amsmath}
\usepackage{amsfonts,amssymb,mathrsfs,amscd}
\usepackage{verbatim}
\usepackage{eucal}
\usepackage{graphicx}
\usepackage{enumerate}
\usepackage[hyper]{msb-a}

\JournalName{%%Математический сборник
}

\numberwithin{equation}{section}
\theoremstyle{plain}
\newtheorem{theorem}{Теорема}
\newtheorem{maintheorem}{Основная теорема}
\newtheorem{propos}{Предложение}%%[subsection]
\newtheorem{ThB}{Теорема B}
\newtheorem{ThA}{Теорема A}
\newtheorem{ThL}{Теорема Лиувилля}

\theoremstyle{definition}

\newtheorem{proof}{Доказательство}
\newtheorem{remark}{Замечание}%%[subsection]
\renewcommand{\leq}{\leqslant} 
\renewcommand{\geq}{\geqslant}

\newcommand{\RR}{\mathbb{R}} 
\newcommand{\CC}{\mathbb{C}} 
\newcommand{\NN}{\mathbb{N}} 
\newcommand{\oB}{{\overline B}}

\DeclareMathOperator{\ord}{{\sf ord}}

\DeclareMathOperator{\dd}{\,{\mathrm d\!}}

\DeclareMathOperator{\sbh}{{\sf sbh}}

\begin{document} 
\title{Теоремы типа Лиувилля с ограничениями вне малых множеств на окружностях или сферах для функций конечного порядка}
	
\author[B.\,N.~Khabibullin]{Б.\,Н.~Хабибуллин}
\address{Башкирский государственный университет}
\email{khabib-bulat@mail.ru}

\date{04.09.2020}
\udk{517.574 : 517.576 : 517.550.4 : 517.547.2 : 517.518.244}

 \maketitle

\begin{fulltext}

\begin{abstract}  Доказано, что субгармонические функции конечного порядка на конечномерном вещественном  пространстве, ограниченные сверху вне некоторых асимптотически малых множеств на сферах, ограничены сверху всюду. Отсюда следует, что субгармонические функции конечного порядка на комплексной плоскости, целые и плюрисубгармонические функции конечного порядка, а также выпуклые или гармонические функции конечного порядка, ограниченные сверху вне таких же множеств на сферах, являются постоянными. Наши результаты и методы доказательства новые и для функций одной комплексной переменной.

Библиография:   12 названий 

Ключевые слова: целая функция, (плюри)субгармоническая функция, выпуклая функция, среднее по окружности и сфере, функция конечного порядка, теорема Лиувилля  

We prove that subharmonic functions of finite order on finite dimensional real space, bounded from above outside of some asymptotically small sets on spheres, are bounded from above everywhere. It follows that subharmonic functions  of finite order  on the complex plane, entire and  plurisubharmonic functions  of finite order, and convex or harmonic  functions  of finite order  bounded from above outside of such sets  on spheres are constant. Our results and methods of proof are also new for functions of one complex variable.

Bibliography:  12 titles 



\end{abstract}

\markright{Теоремы типа Лиувилля с ограничениями вне малых множеств\dots}

\footnotetext[0]{Исследование выполнено за счёт гранта Российского научного фонда (проект № 18-11-00002).}
	
\section{Введение}
\subsection{Предшествующие результаты}\label{s10}
Всюду в статье  $\NN:=\{1,2,\dots\}$ --- множество {\it натуральных чисел,\/}  $\RR$ и $\CC$ ---  соответственно поля {\it вещественных\/} и {\it комплексных чисел.\/} Для чисел  $m, n\in \NN$ 
векторные пространства $\RR^m$ над $\RR$ и $\CC^n$ над $\CC$ рассматриваются как евклидовы пространства 
с евклидовой нормой-модулем  $|\cdot|$. При необходимости и возможности пространство  $\CC^n$ отождествляется с 
$\RR^{2n}$. Исходный классический результат --- 
 
\begin{ThL}[\cite{Rans}, \cite{HK}, \cite{HorC}, \cite{ABR}] Из ограниченности сверху выпуклой или гармонической функции на $\RR^m$, а также целой  или плюрисубгармонической функции на $\CC^n$ следует, что эта функция постоянна.   
\end{ThL}
%%Такое же заключение верно и для {\it ограниченных сверху\/} субгармонических функций на $\CC$ \cite[следствие %%2.3.4]{Rans} и, как очевидное  следствие,  плюрисубгармонических функций на $\CC^n$,  выпуклых функций на {\it %%вещественной прямой\/} $\RR$ и, как мгновенное следствие,  на $\RR^m$ при $1<m\in \NN$,  гармонических %%функций  на $\RR^m$ при любых $m\in \NN$ \cite[теорема 1.19]{HK}.

При условии ограниченности сверху целой функции одной комплексной переменной  вне малых множеств    
известна и находит применения следующая 
\begin{ThB}[{\rm (\cite[лемма 4.2]{BarBelBor18},  \cite[лемма 4.2]{BarBelBor18a},  \cite[теорема 2.1]{AlemBarBelHed20})}]  Если  целая функция   конечного порядка  на $\CC$ ограничена вне  
 $E\subset \CC$, а  %%множество $E$  нулевой плоской плотности  в том смысле, что
площади пересечений $E$ с кругами радиуса $r\in \RR^+:=\{x\in \RR\colon x\geq 0\}$ с центром  в нуле определены и являются величиной порядка $o(r^2)$ при $r\to +\infty$, то эта функция  постоянная.
\end{ThB}

<<Быстрое>> доказательство теоремы B изложено в \cite{Kha20ar}. В \cite{Kha20arU}--\cite{Kha20U} 
теорема B перенесена на функции многих переменных:

\begin{ThA}[{\rm \cite[теорема 1]{Kha20arU}, \cite[теорема 1]{Kha20U}}]  
Если  плюрисубгармоническая или целая функция  конечного порядка на $\CC^n$ ограничена  сверху
вне  $E\subset \CC^n$, а объёмы  пересечений $E\subset \CC^n$ с шарами радиуса $r\in \RR^+$ с центром  в нуле определены и являются величиной порядка $o(r^{2n})$ при $r\to +\infty$, то  эта функция  постоянная.
\end{ThA}

Из  теоремы \ref{th1} настоящей статьи следует, что заключение теоремы B сохраняется и после замены ограничения  на множество $E$ из теоремы B на существенно более слабое: 
\begin{enumerate}
\item[{[{\bf E}]}] {\it существует неограниченная возрастающая\footnote{{\it Положительность\/}  всегда понимается как $\geq 0$, а {\it возрастание\/} определяется с нестрогими неравенствами. То же самое по  {\it отрицательности,\/}
{\it убыванию,\/} {\it монотонности}. Со строгими неравенствами --- это {\it строгая положительность,\/}
{\it строгое возрастание\/} и т.д., и т.п.} последовательность  положительных чисел $(r_k)_{k\in \NN}$
на $\RR^+$, удовлетворяющая условию 
\begin{equation}\label{rk}
\limsup_{k\to \infty}\frac{r_{k+1}}{r_k}<+\infty,
\end{equation} 
 для которой  множества  $e_k:=\{\arg z\colon z\in E, |z|= r_k\}\bigcap  [0,2\pi)$ измеримы 
по  линейные мере Лебега на интервале $[0,2\pi)$, а предел  линейных мер Лебега множеств $e_k$ при $k\to \infty$ равен нулю.} 
\end{enumerate}

Более того, методы настоящей работы дают  подобные результаты и для  плюрисубгармонических и целых  функций на $\CC^n$ при всех  $n\in \NN$, а также  для выпуклых и гармонических функций на $\RR^m$ при всех $m\in \NN$.  Далее при обращении к (плюри)субгармоническим, гармоническим,  целым и выпуклым  функциям и их свойствам  вполне достаточны основные сведения о них из \cite{Rans}, \cite{HK}, \cite{HorC}, \cite{ABR}.

\subsection{Средние по сферам и шарам. Порядок функции}\label{IAmf}
 Через   
\begin{subequations}\label{SB}
\begin{align}
 B_m(x,r)&:=\bigl\{x' \in \RR^m \colon |x'-x|< r\bigr\},
\tag{\ref{SB}B}\label{SBB}\\
\overline  B_m(x,r)&:=\bigl\{x' \in \RR^m \colon |x'-x|\leq r\bigr\},
\tag{\ref{SB}$\rm \overline B$}\label{SBoB}\\
 S_{m-1}(x,r)&:=\bigl\{x' \in \RR^m \colon |x'-x|= r\bigr\}=\overline  B_m(x,r)\!\setminus\!B_m(x,r),
\tag{\ref{SB}S}\label{SBS}
\end{align}
\end{subequations}
обозначаем соответственно {\it открытый\/} и  {\it замкнутый шары,\/} а также  {\it  сферу\/}  в $\RR^m$ {\it радиуса $r\in \RR^+$ с центром $x\in \RR^m$,\/}, а $B_m(r):=B_m(0,r)$, $\oB_m(r):=\oB_m(0,r)$ и $S_{m-1}(r):=S_{m-1}(0,r)$, где нижние индексы 
$m$ и $m-1$ указывают на {\it размерность\/.} Такие же обозначения  с $m:=2n$ используем  в  $\CC^n$, отождествляемом с $\RR^{2n}$. 
Нижние индексы $m$ и $m-1$ по возможности опускаем, когда это не вызывает разночтений. 
{\it Меру Лебега\/} на $\RR^m$ %%и её сужения на подмножества в шары $\overline  B(x,r)$ и $\overline  B(r)$
 обозначаем через $\lambda $, а {\it поверхностную меру\/} на $S(r)$  через $\sigma_r$.
Через $b_m$ и $s_{m-1}$ обозначаем соответственно {\it объём\/} $\lambda\bigl(B(1)\bigr)$ {\it единичного шара\/} $\oB(1)$ в $\RR^m$ и {\it площадь\/ $\sigma_r\bigl(S(r) \bigr)$ единичной сферы\/} $S(1)\subset \oB(1)\subset \RR^m$.  
Для  $\lambda$-интегрируемой функции $v\colon \overline B_m(x,r)\to  \overline \RR$   её {\it среднее по шару\/}
$\overline B_m(x,r)$ обозначаем как 
\begin{subequations}\label{BS}
\begin{align}
{\bf B }_v(x,r)&:=
%%\frac{1}{\lambda\bigl(\overline B(x,r)\bigr)}\int_{\overline B(x,r)} v \dd \lambda=
\frac{1}{b_mr^m}\int_{\overline B(x,r)} v \dd \lambda,
&{\bf B }_v(r):={\bf B }_v(0, r),
\tag{\ref{BS}B}\label{{BS}B}
\\
\intertext{а для  $\sigma_r$-интегрируемой $v\colon S_{m-1}(x,r)\to  \overline \RR$ её {\it среднее по сфере\/}
$S_{m-1}(x,r)$ как}
{\bf S }_v(x,r)&:=
%%\frac{1}{\sigma_r\bigl(S(x,r)\bigr)}\int_{S(x,r)} v \dd \sigma_r=
\frac{1}{s_{m-1}r^{m-1}}
\int_{S(r)} v(x+y) \dd \sigma_r(y), 
&{\bf S }_v(r):={\bf S }_v(0, r).
\tag{\ref{BS}S}\label{{BS}S}
\end{align}
\end{subequations}
Для  $\lambda$-интегрируемой функции $v\colon \overline B(x,r)\to  \overline \RR$ средние по сферам  \ref{{BS}S}
определены почти всюду по линейной мере Лебега на $[0,r]$ и связаны  равенством %%\cite[(1)]{Br}
\begin{equation*}
{\bf B }_v(x,r)=\frac{m}{r^m}\int_0^r{\bf S}_v(x,t)t^{m-1}\dd t.
\end{equation*}

Через $\sbh(S)$ обозначаем класс всех {\it субгармонических\/} ({\it локально выпуклых} при $m=1$) {\it функций\/} на каких-либо открытых  окрестностях множества $S\subset \RR^m$.   
Роль средних из \eqref{BS}  для субгармонических функций $v\not\equiv -\infty$ обусловлена полностью характеризующими их, при условии полунепрерывности сверху  и локальной интегрируемости по  $\lambda$, неравенствами  для  средних по шару и сфере \cite{Rans}, \cite{HK}. В частности,  
с постоянной 
\begin{equation*}
a_m:=\begin{cases}
1/2 &\text{при $m=1$},\\
1/\sqrt{e} &\text{при $m=2$},\\
1/\sqrt[{m-2}]{m/2} &\text{при $m\geq 3$}
\end{cases}
\end{equation*}
имеют место точные неравенства \cite[теорема]{Br}, \cite[(5)--(6), перед теоремой 3]{FM}
\begin{equation}\label{inav}
v(0)\leq  {\bf S}_v(a_mR)\leq {\bf B}_v(R)\leq {\bf S }_v(R) %%\quad\text{для  всех $v\in \sbh\bigl(\oB(x,r)\bigr)$}.
\end{equation}
{\it для любой функции  $v\colon \overline B(R)\to  \overline \RR$ c четырьмя свойствами   
\begin{enumerate}[{\rm [$1v$]}]
\item\label{1v} функция $v$ интегрируема на $\overline B(R)$ по мере Лебега $\lambda$;
\item\label{2v} сужение $v\bigm|_{B(R)}$ функции $v$ на $B(R)\subset \RR^m$
--- субгармоническая функция;
\item\label{3v} 
сужение $v\bigm|_{S(R)}$ функции $v$ на $S(R)\subset \RR^m$ интегрируемо по  мере $\sigma_R$;
\item\label{4v} 
$\limsup\limits_{B(R)\ni x'\to x} v(x')\leq v(x)$ для $\sigma_R$-почти каждой точки $x\in S(R)$.
\end{enumerate}
}

Пусть функция $M$ со значениями в  {\it расширенной вещественной прямой\/} $\overline \RR:=\RR\cup\{\pm\infty\}$ определена   на  $\RR^+$, в $\RR^m$ или в $\CC^n$, но, возможно,   вне некоторого шара $B(r)$ радиуса $r\in \RR^+$.   {\it Порядок\/}  функции $M$ (около $\infty$)  определяется  как
\cite[2.1]{KhI} 
\begin{equation}\label{ord}
\ord[f]:=\limsup_{|x|\to +\infty}\frac{\ln \bigl(1+M^+(x)\bigr)}{\ln |x|}\in \RR^+\cup \{+\infty\}, 
\end{equation} 
где $M^+\colon x\mapsto \max\bigl\{0, M(x)\bigr\}$ --- {\it положительная часть функции} $M$. Порядок целой функции $f$ на $\CC^n$ или голоморфной функции $f$ вне некоторого шара  --- это порядок $\ord\bigl[\ln |f|\bigr]$ {\it плюрисубгармонической\/} функции $\ln|f|$ при соглашении, что субгармонические функции на подмножествах комплексной плоскости  $\CC$ также называем и плюрисубгармоническими. 

\subsection{Основные результаты} 

\begin{maintheorem}
Пусть $m\in \NN$,  $r_0\in \RR^+$, $E\subset \RR^m$,   $(r_k)_{k\in \NN}$  --- 
последовательность из\/ {\rm [{\bf E}]}
с $r_1>r_0$ и ограничением  \eqref{rk}, 
для которой  множества $E_k:=E\cap S(r_k)\subset S(r_k)$ измеримы по поверхностной мере $\sigma_{r_k}$  и 
\begin{equation}\label{Ek}
\limsup_{k\to \infty} \frac{\sigma_{r_k}(E_k)}{r_k^{m-1}}=0. 
\end{equation}
Тогда для любой функции $v\in \sbh \bigl(\RR^m\!\setminus\!\oB (r_0)\bigr)$ конечного порядка $\ord[v]<+\infty$ имеет место равенство
\begin{equation}\label{supk}
%%\sup_{\RR^m\!\setminus\!\oB (r_0)} v=\sup\biggl\{\sup_{\oB(r_1)\!\setminus\!\oB(r_0)}v, \sup_{\RR^m\!\setminus\!E} v \biggr\}. 
\sup_{\RR^m\!\setminus\!\oB (r_0)} v:=\max\{M_E, M_0\}, 
\quad\text{где $M_E:=\sup_{\RR^m\!\setminus\!E} v$, \quad  $M_0:=\sup_{\oB(r_1)\!\setminus\!\oB(r_0)}v$.} 
\end{equation} 
\end{maintheorem}
\begin{remark}\label{rem1} Чем больше исключаемое множество $E\subset \RR^m$ в \eqref{supk}, тем заключительное равенство \eqref{supk} сильнее.  Поэтому среди всех множеств $E$ с фиксированными пересечениями  $E_k:=E\cap S(r_k)\subset S(r_k)$  выбор $E$ оптимален, если включить в $E$  все   открытые шаровые слои $B(r_{k+1})\!\setminus\! \oB(r_{k}) $ между сферами  $S(r_k)$ и $S(r_{k+1})$, что никак не может нарушить условие  \eqref{Ek}. 
\end{remark}
\begin{theorem}\label{th1} В  условиях основной теоремы для $\CC^n$, отождествлённого с $\RR^{2n}$ и  $m:=2n$,   ограниченность сверху  на $\oB(r_1)\!\setminus\!\oB(r_0)$ и одновременно на  $\CC^n\!\setminus\!E$   плюрисубгармонической или голоморфной функции на $\CC^n\!\setminus\!\oB (r_0)$ конечного порядка влечёт за собой  её постоянство всюду   на $\CC^n\!\setminus\!\oB (r_0)$.
\end{theorem}
\begin{remark}\label{rem2}
По  теореме \ref{th1} действительно получаем заключение  теоремы B при более слабом ограничении [{\bf E}], обсуждавшееся выше  после её формулировки, поскольку  обозначения из  [{\bf E}] и \eqref{Ek}
согласованы  при $n=1$ и $m=2\cdot 1$ для  меры Лебега $\lambda$ на $[0,2\pi)$
и меры длины окружности $\sigma_{r_k}$  на $S_1(r_k)$ как 
\begin{equation}\label{x}
\lambda(e_k)=\frac{\sigma_{r_k}(E_k)}{r_k}, 
\quad  \limsup_{k\to \infty} \lambda(e_k)=\limsup_{k\to \infty} \frac{\sigma_{r_k}(E_k)}{r_k}. 
\end{equation}
\end{remark}

\begin{remark}
При $n>1$ по известному принципу Хартогса  функция, голоморфная на дополнении $\CC^n\!\setminus\!\oB (r_0)$ компакта $\oB (r_0)$, является  сужением некоторой  целой функции на $\CC^n\!\setminus\!\oB (r_0)$. %%Следовательно, в этом случае получаем  её 
\end{remark}

\begin{theorem}\label{th2} Пусть для каждой точки ${\sf s}\in S(1)$ и соответствующей ей комплексной прямой $\CC_{\sf s}:=\{z{\sf s}\colon z\in \CC\}$ в $\CC^n$,  рассматриваемой как комплексная плоскость, найдутся  $E_{\sf s}\subset \CC_{\sf s}$ и последовательность $\bigl(r_k({\sf s})\bigr)_{k\in \NN}$, для которых выполнено\/ {\rm [{\bf E}]}. 
Если для плюрисубгармонической или целой функции на $\CC^n$
её сужение на каждую плоскость $\CC_{\sf s}$ --- функция конечного порядка, ограниченная сверху на  $\CC_{\sf s}\!\setminus\!E_{\sf s}$ своей постоянной, то эта функция  постоянна на всём\/ $\CC^n$. 
\end{theorem}

\begin{theorem}\label{th3} В  условиях основной теоремы  ограниченность сверху   на дополнении $\RR^m\!\setminus\!E$  выпуклой или гармонической функции на $\RR^m$ конечного порядка  влечёт за собой  её постоянство всюду на $\RR^m$.
\end{theorem}

\section{Неравенства со средними по сферам и шарам}
%% для субгармонических  функций}
Установим вспомогательные неравенства в несколько большем объёме, чем необходимо в настоящей статье, в частности, с целью сравнить подходы к оценкам субгармонических функций  через средние по сферам и средние по шарам, а также создать задел для исследования субгармонических функций и их подклассов в круге и шаре.  
\begin{propos}\label{1} Пусть $0<r<R\in \RR^+$,  $v$ ---  функция на   $\overline B(R)\subset \RR^m$ с 
четырьмя свойствами {\rm [$1v$]--[$4v$]} %%{inav}
и положительной частью $v^+$.
%%, субгармоническая в открытом шаре $B(R)\subset \RR^m$, 
%%сужение $v\bigm|_{S(r)}$ функции $v$ на $S(r)$ интегрируемо по поверхностной мере $\sigma_R$ и 
%%\begin{equation*}
%%\limsup_{B(R)\ni x'\to x} v(x')\leq v(x)\quad\text{для каждой точки $x\in S(R)$}.
%%\end{equation*}  
Тогда 
\begin{subequations}\label{BB}
\begin{align}
v(x)&\overset{\eqref{inav}}{\leq} {\bf  S}_{v}\bigl(x,a_m(R-r)\bigr) \leq {\bf  S}_{v^+}\bigl(x,a_m(R-r)\bigr)
\overset{\eqref{inav}}{\leq} {\bf  B}_{v^+}(x,R-r)
\tag{\ref{BB}B}\label{{BB}B}
\\
&\leq \Bigl(1+\frac{r}{R-r}\Bigr)^m  {\bf B}_{v^+}(R) 
\overset{\eqref{inav}}{\leq} \Bigl(1+\frac{r}{R-r}\Bigr)^m  {\bf S}_{v^+}(R) 
\text{ при   $x\in \oB(r)$};
\tag{\ref{BB}S}\label{{BB}S}\\
v(x)&\overset{\eqref{inav}}{\leq}  {\bf S}_v(x,t)\leq 
\frac{\bigl(R+r+t\bigr)R^{m-2}}{\bigl(R-r-t\bigr)^{m-1}}
{\bf S}_{v^+}(R)
\leq\frac{4R^{m-1}}{\bigl(R-r-t\bigr)^{m-1}}
{\bf S}_{v^+}(R)
\tag{\ref{BB}t}\label{{BB}t}
\\
&=4 \Bigl(1+\frac{r+t}{R-(r+t)}\Bigr)^{m-1} {\bf S}_{v^+}(R) \quad\text{при $x\in \oB(r)$ и\/ $0<t<R-r$,} 
\tag{\ref{BB}r}\label{{BB}r}\\
\intertext{откуда при выборе $t:=\frac12(R-r)$ %%из  \eqref{{BB}t}--\eqref{{BB}r} 
имеем}
v(x)&
\overset{\eqref{inav}}{\leq} {\bf S}_v\Bigl(x,\frac{R-r}{2}\Bigr)
\leq  2^{m+1}\Bigl(1+\frac{r}{R-r}\Bigr)^{m-1} {\bf S}_{v^+}(R)
 \quad\text{при $x\in \oB(r)$}, 
\tag{\ref{BB}R}\label{{BB}R}\\
\intertext{а  \eqref{{BB}r} при  $0<t\to 0$ даёт неравенство}
v(x)&
\leq  4\Bigl(1+\frac{r}{R-r}\Bigr)^{m-1} {\bf S}_{v^+}(R)
 \quad\text{при $x\in \oB(r)$}. 
\tag{\ref{BB}v}\label{{BB}v}
\end{align}
\end{subequations}
\end{propos}

\begin{proof} Ссылки-метки над тремя неравенствами в \eqref{BB} обосновывают их. Второе  неравенство в 
\eqref{{BB}B}  --- очевидно следствие положительности $v^+\geq v$. Переходное неравенство от 
 \eqref{{BB}B} к  \eqref{{BB}B} установлено в \cite[предложение 1]{Kha20arU}--\cite[предложение 1]{Kha20U} совершенно элементарным приёмом для произвольных положительных $\lambda$-измеримых  функций $v^+$ на $\oB(r)$.

Прежде чем доказывать  неравенства   \eqref{{BB}t}--\eqref{{BB}r} для корректности построений и выкладок  отметим, что, не умаляя общности, можем рассмотреть только  функции $v\in \sbh \bigl(\oB(R)\bigr)$, к которым можно перейти с помощью гомотетии переменной с коэффициентом $<1$, а затем вернуться к исходной функции $v$, устремляя этот 
коэффициент в   \eqref{{BB}t}--\eqref{{BB}r} к единице. Ввиду непрерывности средних по шару и сфере для субгармонических функций неравенства  \eqref{{BB}t}--\eqref{{BB}r} при таком переходе не нарушатся. Более того, достаточно рассматривать только {\it непрерывные положительные  функции\/} $v\in \in \sbh \bigl(\oB(R)\bigr)$, а потом
для произвольной функции $v^+\in \in \sbh \bigl(\oB(R)\bigr)$ воспользоваться убывающей к ней последовательностью таких функций.  Для положительной непрерывной функции $v^+\in \sbh \bigl(\oB(R)\bigr)$ уже можно рассмотреть её {\it положительное гармоническое продолжение\/} $H\geq v^+$ со сферы $S(R)$ внутрь $\oB(R)$, т.\,е. в $B(R)$,  построенное с помощью интеграла Пуассона. Для такой функции $H$ по неравенству Харнака 
\cite[теорема 1.18]{HK} имеют место неравенств 
\begin{equation*}
 H (x')\leq \frac{(R+|x'|)R^{m-2}}{(R-|x'|)^{m-1}}H(0)
\quad\text{для любой точки $x'\in B(R)$}
\end{equation*}
где $H(0)$ равно  среднему ${\bf S}_{v^+}(R)$ по сфере $S(R)$ функции $v^+$, а $v^+(x')\leq H(x')$, поскольку $H$ --- мажоранта 
функции $v^+$. Таким образом, 
\begin{equation*}
v^+(x')\leq \frac{(R+|x'|)R^{m-2}}{(R-|x'|)^{m-1}}{\bf S}_{v^+}(R)
\quad\text{для любой точки $x'\in B(R)$}.
\end{equation*}
Дробь в правой части здесь возрастает при росте $|x'|$. Следовательно, полагая $x'=x+y$
для $x\in \oB(r)$ и $y\in S(t)$, при фиксированных $0<r<R$ и $0<t<R-r$, имеем  $|x'|\leq r+t<R$ и   
\begin{equation*}
v^+(x+y)\leq \frac{\bigl(R+(r+t)\bigr)R^{m-2}}{\bigl(R-(r+t)\bigr)^{m-1}}{\bf S}_{v^+}(R)
\quad\text{при  $x\in \oB(r)$ и $y\in S(t)$}.
\end{equation*}
Правая часть неравенство здесь --- это фиксированное число, и интегрирование обеих частей этого неравенства по вероятностной мере 
$\frac{1}{s_{m-1}t^{m-1}}\sigma_t$ по определению \eqref{{BS}S} среднего по сфере даёт  
\begin{equation*}
{\bf S}(x,t)=\frac{1}{s_{m-1}t^{m-1}}\int_{S(t)}v^+(x+y)\dd \sigma_t(y)\leq \frac{(R+r+t)R^{m-2}}{(R-r-t)^{m-1}}{\bf S}_{v^+}(R)
\end{equation*}
при  всех $x\in \oB(r)$ и $0<t<R-r$, что доказывает второе неравенство в \eqref{{BB}t}. 
Продолжения этой цепочки неравенств от  \eqref{{BB}t} к \eqref{{BB}r} и до \eqref{{BB}R} элементарно.
\end{proof}

\begin{remark}\label{rem4} Простой сравнительный анализ неравенств \eqref{{BB}B}--\eqref{{BB}S} с неравенствами \eqref{{BB}t}--\eqref{{BB}R}  показывает, что первые могут быть более эффективными при 
 $0<r\to 0$ за счёт стремления множителя в них перед ${\bf S}_{v^+}(R)$ к единице, в то время как во второй группе неравенств этот множитель всегда не меньше четырёх. Напротив, при $R>r\to R$ ситуация противоположная, поскольку 
порядок роста множителя перед ${\bf S}_{v^+}(R)$ в \eqref{{BB}B}--\eqref{{BB}S}  --- это  $O\bigl((R-r)^{-m}\bigr)$  при $R-r\to 0$, а в остальных неравенствах --- только   $O\bigl((R-r)^{-m+1}\bigr)$.
В соизмеримых с $R$ границах $q\leq R/r\leq Q$ со строго положительными  числами   $q>1$ и $Q<+\infty$ пара этих групп неравенств равнозначна, если отвлечься от размерности $m$, и даёт постоянную, зависящую  только от $q$ и $Q$. В настоящей статье используется лишь последний <<соизмеримый>> вариант. А вариант с малыми $r$ представляется более приспособленным к случаю медленно растущих исследуемых функций, в то время как вариант с $r$, близкими к  $R$ ---  для исследования быстро растущих функций в $\RR^m$ или в $\CC^n$, а также  для применений к функциям в круге  или в шаре.      
\end{remark}

\begin{propos}\label{2}
 Пусть $0<r<R\in \RR^+$,  $v$ ---  функция на   $\overline B(R)\subset \RR^m$ с 
четырьмя свойствами {\rm [$1v$]--[$4v$]} %%{inav}
и положительной частью $v^+$, а подмножество $E\subset S(r)\subset \RR^m$ 
измеримо по поверхностной мере $\sigma_r$ на $S(r)$. Тогда 
\begin{equation}\label{v+}
%%{\sf B}_v(r)\leq \frac{1}{b_mr^m}\int_{B(r)\!\setminus\!E}v\dd \lambda
\int_{E}v\dd \sigma_r \leq 
\min\Bigl\{4, 
1+\frac{r}{R-r}\Bigr\}
\Bigl(1+\frac{R+r}{R-r}\Bigr)^{m-1} 
\sigma_r (E){\bf S}_{v^+}(R).
%%\int_{\overline B(R)}v^+\dd \lambda
\end{equation}
\end{propos}
\begin{proof} %%В \cite[предложение 2]{Kha20ar} 
Интегрирование крайних частей неравенств \eqref{{BB}B}--\eqref{{BB}S}  по поверхностной  мере  $\sigma_r$ и по  множеству $E\subset S(r)\subset \oB(r)$ даёт %%неравенство 
\begin{equation*}
\int_{E}v\dd \sigma_r \leq 
\int_E\Bigl(1+\frac{r}{R-r}\Bigr)^{m}
{\bf S}_{v^+}(R) \dd \sigma_r
=\Bigl(1+\frac{r}{R-r}\Bigr)^{m}
{\bf S}_{v^+}(R) \sigma_r (E), 
\end{equation*}
а такая  же  процедура применительно к неравенству \eqref{{BB}v} ---
\begin{equation*}
\int_{E}v\dd \sigma_r \leq 
\int_E 4\Bigl(1+\frac{r}{R-r}\Bigr)^{m-1} {\bf S}_{v^+}(R) \dd \sigma_r
= 4\Bigl(1+\frac{r}{R-r}\Bigr)^{m-1} {\bf S}_{v^+}(R)
\sigma_r (E). 
\end{equation*}
Пересечение  этих соотношение --- в точности \eqref{v+}. 
\end{proof}
\begin{remark} Прежде всего и к \eqref{v+} всецело относится предыдущее замечание \ref{rem4}. 
Кроме того, при большой размерности $m$ преимущество у второй компоненты в минимуме
из  \eqref{v+}, но часто может быть наоборот при малых $m$.
\end{remark}

\section{Доказательства основных результатов}
\begin{proof}[основной теоремы] Если хотя бы одна из величин $M_0$ или $M_E$ равна $+\infty$, то 
\eqref{supk} очевидно. В противном случае положим  $M:=\max\{M_0,M_E\}\in \RR$, и сначала 
рассмотрим случай  функции $v$, субгармонической на всём пространстве $\RR^m$, для которой по принципу максимума  
\begin{equation}
\sup_{\oB(r_1)} v\overset{\eqref{supk}}{\leq} M_0\leq M, \quad\text{а также $\sup_{S(r_k)\!\setminus\!E_k} v\leq M$} 
\end{equation}
По  условию  \eqref{rk} нетрудно подобрать числа  $q>1$ и $Q\in \RR^+$   и выделить такую строго возрастающую подпоследовательность из последовательности $(r_k)_{k\in \NN}$, что отношение каждого последующего члена подпоследовательности к предшествующему попадает в отрезок $[q, Q]$, а начинается она по-прежнему с $r_1$.
 Сохраним за этой подпоследовательностью то же обозначение   $(r_k)_{k\in \NN}$, для которой теперь 
\begin{equation}\label{r}
q\leq \frac{r_{k+1}}{r_k}\leq Q\quad\text{для каждого $k\in \NN$},
\end{equation}
а для соответствующих множеств $E_k=E\cap D(r_k)$ по-прежнему выполнено \eqref{Ek}.
Вместо функции $v\in \sbh(\RR^m)$ рассмотрим положительную  функцию
\begin{equation}\label{v-M+}
V:=(v-M)^+=V^+\in \sbh(\RR^m)
\end{equation} 
 --- положительную часть функции $v-M$, которая {\it отрицательна на\/} $S(r_1)$ и {\it на всех множествах\/} 
$S(r_k)\!\setminus\!E_k$ при $k\in \NN$, а также конечного порядка.  Интегрирование по поверхностной мере $\sigma_{r_k}$ для каждого $k\in \NN$
даёт неравенство 
\begin{equation*}\label{Sv+}
{\bf S}_{V}(r_k)=\frac{1}{s_{m-1}r_k^{m-1}} \int_{S(r_k)}V\dd \sigma_{r_k}\leq 
\frac{1}{s_{m-1}r_k^{m-1}} \int_{E_k}V\dd \sigma_{r_k}.
\end{equation*}
Применим теперь к интегралу в правой части предложение \ref{2} с $r_k$ в роли  $r$, $r_{k+1}$ в роли $R$,  $E_k$ в роли $E$ и функцией $V=V^+$ в роли $v$:
\begin{multline*}
%%\frac{1}{s_{m-1}r_k^{m-1}}
 \int_{E_k}V\dd \sigma_{r_k}\leq 
\min\Bigl\{4, 1+\frac{r_k}{r_{k+1}-r_k}\Bigl\}\Bigl(1+\frac{r_k}{r_{k+1}-r_k}\Bigr)^{m-1}\sigma_{r_k}(E_k){\bf S}_{V}(r_{k+1})\\
\overset{\eqref{r}}{\leq}  4\Bigl(\frac{q}{q-1}\Bigr)^{m-1}
\sigma_{r_k}(E_k){\bf S}_{V}(r_{k+1}).
\end{multline*}
Применяя эту оценку к предшествующей, получаем  
\begin{equation*}
{\bf S}_V(r_k)
{\leq} \frac{4}{s_{m-1}} \Bigl(\frac{q}{q-1}\Bigr)^{m-1}
\frac{\sigma_{r_k}(E_k)}{r_k^{m-1}}{\bf S}_V(r_{k+1}),
\end{equation*}
где ${\bf S}_V(r_{k+1})\leq {\bf S}_V(Qr_{k})$ ввиду \eqref{r} и возрастания среднего по сферам ${\bf S}_V$.  
Отсюда согласно условию \eqref{rk}  имеем
\begin{equation*}
\lim_{k\to +\infty}\frac{{\bf S}_V(r_k)}{{\bf S}_V(Qr_k)} =0,
\end{equation*}
где функция  ${\bf S}_V$ конечного порядка, поскольку таковой является функция  $V$. Но последнее предельное равенство нулю  для положительной возрастающей функции конечного порядка на $\RR^+$ на последовательности
$(r_k)_{k\in \NN}$, удовлетворяющей \eqref{rk}, возможно лишь для нулевой функции   ${\bf S}_V\equiv 0$ на $\RR^+$. 
Таким образом, $(v-M)^+\equiv 0$ и  $v\leq M=\max\{M_0,M_E\}$ . Отсюда ввиду определения чисел $M_0$ и $M_E$
в \eqref{supk} получаем требуемое первое равенство в \eqref{supk}. 

В случае субгармонической функции $v$, определённой  лишь на $\RR^m\!\setminus\!\oB(r_0)$, можем заменить  её на функцию 
\begin{equation}\label{sM}
\begin{cases}
M_0 &\text{на $\oB(r_1)$},\\
\sup\{v,M_0\} &\text{на $\RR_m\!\setminus\!\oB(r_1)$},
\end{cases}
\end{equation}
которая {\it субгармонична\/} уже {\it на всём\/} $\RR^m$ и обладает всеми требуемыми в основной теореме свойствами, но с числом $\max\{M_0,M_E\}$ в роли  $M_E$, что не меняет заключения  основной теоремы.
\end{proof}

\begin{proof}[теоремы \ref{th1}] 
Рассмотрим  плюрисубгармонические функции $v$ на $\CC^n\!\setminus\!\oB (r_0)$. При этом можно перейти к функции \eqref{sM}, которая {\it плюрисубгармонична\/} на всём $\CC^n$, и, тем более субгармонична
на $\RR^{2n}$, отождествлённом с $\CC^n$.  
При этом выполнены все условия основной теоремы, следовательно функция \eqref{sM} ограничена сверху всюду на $\CC^n$. Но поскольку она плюрисубгармонична, то по теореме Лиувилля для плюрисубгармонических функций, получаем, что она постоянна.

Если $f$ --- голоморфная функция на    $\CC^n\!\setminus\!\oB (r_0)$ из теоремы \ref{th1}, то для плюрисубгармонической функции $\ln|f|$ выполнены все условия теоремы \ref{th1}. Следовательно,  $|f|$ --- постоянная функция, откуда и функция $f$ постоянна.
\end{proof}

\begin{proof}[теоремы \ref{th2}] По теореме \ref{th1} для субгармонических функций на комплексной плоскости сразу получаем, что сужение исходной плюрисубгармонической функции на каждую комплексную плоскость 
$\CC_{\vec{e}}$ --- постоянная функция. Но  все комплексные прямые  $\CC_{\vec{e}}$  имеют общую точку $0$, следовательно, и исходная плюрисубгармонической функция постоянна.
От  исходной целой функции $f$  переходим к плюрисубгармонической функции $\ln|f|$, которая по доказанной части 
теоремы \ref{th2} постоянна, откуда постоянна и $f$.
\end{proof}

\begin{proof}[теоремы \ref{th3}] Выпуклая или гармоническая функция $v$ субгармонична, а  по основной теореме эта  функция  ограниченна сверху на всём $\RR^m$. По теореме Лиувилля она постоянна. 
\end{proof}

\begin{remark}
При $m=2$ от гармонической функции  в теореме \ref{th3} достаточно требовать её гармоничности не всюду на плоскости, а лишь в проколотой плоскости $\CC\!\setminus\!\{0\}$ \cite[следствие 3.3]{ABR}. Для иных размерностей это не верно.
\end{remark}

\end{fulltext}

\end{document}